
\baselineskip=14pt
\parskip=10pt

\magnification=\magstephalf

\def\1{{\overline{1}}}
\def\2{{\overline{2}}}
\parindent=0pt
\overfullrule=0in

\def\frac#1#2{{#1 \over #2}}
\centerline
{\bf 
On the Intriguing Problem of Counting $(n+1,n+2)$-core partitions into Odd Parts
} 
\rm
\centerline
{\it Anthony ZALESKI and Doron ZEILBERGER}

{\bf Abstract}: Tewodros Amdeberhan and Armin Straub initiated the study of 
enumerating subfamilies of the set of $(s,t)$-core partitions.
While the enumeration of $(n+1,n+2)$-core partitions into {\it distinct} parts is relatively easy
(in fact it equals the Fibonacci number $F_{n+2}$),
the enumeration of $(n+1,n+2)$-core partitions into {\it odd} parts remains elusive.
Straub computed the first eleven terms 
of that sequence,  
and asked for a ``formula," or at least a fast way, to compute many terms.
While we are unable to find a ``fast" algorithm, we did manage to find a ``faster" algorithm, which enabled
us to compute 23 terms of this intriguing sequence. We strongly believe that this sequence has
an algebraic generating function, since a ``sister sequence" (see the article), is  
OEIS sequence A047749 that  does have an algebraic generating function.
One of us (DZ) is pledging a donation of 100 dollars to the OEIS, in honor of the first person to generate
sufficiently many terms to conjecture (and prove non-rigorously) an algebraic equation for the generating
function of this sequence, and another 100 dollars for a rigorous proof of that conjecture.
Finally, we also develop algorithms that find explicit generating functions for other, more tractable,
families of $(n+1,n+2)$-core partitions.

{\bf Added Jan. 24, 2018}: Paul Johnson observed ([J]) and proved that the ``sister-sequence" counts $(n+1,n+2)$-core partitions into {\bf even} parts,
and even more impressively, related the two sequences, that easily implies a fast way to compute the enumerating sequence that
we were interested in, from which one can easily derive a (rather complicated!) algebraic equation satisfied by the generating function.
We would have needed $53$ terms (rather than $23$) to have guessed it.

{\bf Added Feb. 28, 2018}: Paul Johnson has just posted his beautiful article 

{\it Simultaneous cores with restrictions and a question of Zaleski and Zeilberger}, 

{\tt https://arxiv.org/abs/1802.09621} \quad ,

that does much more than  we asked. A donation to the OEIS, of \$200, in honor of Paul Johnson, has been made.

{\bf Supporting Maple Packages and Output}

All the results in this article were obtained by the use of the Maple packages

$\bullet$ {\tt http://www.math.rutgers.edu/\~{}zeilberg/tokhniot/OddArmin.txt} \quad ,

$\bullet$ {\tt http://www.math.rutgers.edu/\~{}zeilberg/tokhniot/core.txt} \quad ,

$\bullet$ {\tt http://www.math.rutgers.edu/\~{}zeilberg/tokhniot/stCorePlus.txt} \quad ,

whose output files, along with links to diagrams, are available from the {\it front} of this article

{\tt http://www.math.rutgers.edu/\~{}zeilberg/mamarim/mamarimhtml/oddarmin.html} \quad .

{\bf $(s,t)$-Core Partitions}

Recall that a {\it partition} is a non-increasing sequence of positive integers $\lambda=(\lambda_1, \dots, \lambda_k)$ with $k \geq 0$,
called its {\it number of parts}; $n:=\lambda_1 + \dots + \lambda_k$ is called its {\it size},
and we say that $\lambda$ is a {\it partition of $n$}.
Also recall that the {\it Ferrers diagram} 
(or equivalently, using empty squares rather than dots, {\it Young diagram}) 
of a partition $\lambda$ is obtained by placing, in a {\it left-justified} way,
$\lambda_i$ dots at the $i$-th row. For example, the Ferrers diagram of the partition $(5,4,2,1,1)$ is
$$
\matrix
{
* & * & *  &*  & * \cr
* & * & *  &*  &  \cr
* & * &   &  &  \cr
* &  &   &  &  \cr
* &  &   &  &  \cr
} \quad .
$$

Recall also that the {\it hook length} of a dot $(i,j)$ in the Ferrers diagram, $1\leq j \leq \lambda_i$, is the
number of dots to its right (in the same row) plus the number of dots below it (in the same column) plus one
(for itself), in other words $\lambda_i -i+\lambda'_j-j+1$, where $\lambda'$ is the {\it conjugate partition}, 
obtained by reversing the roles of rows and columns. (For example
if $\lambda=(5,4,2,1,1)$ as above, then $\lambda'=(5,3,2,2,1)$.)

Here is a table of hook-lengths of the above partition, $(5,4,2,1,1)$:

$$
\matrix
{
9 & 6 & 4  & 3  & 1 \cr
7 & 4 & 1  &1  &  \cr
4 & 1 &   &  &  \cr
2 &  &   &  &  \cr
1 &  &   &  &  \cr
}  \quad .
$$
It follows that its set of hook-lengths is $\{1,2,3,4,6,7,9\}$. A partition is called an $s$-core if none of its hook-lengths is $s$.
For example, the above partition, $(5,4,2,1,1)$, is a $5$-core, and an $i$-core for all $i \geq 10$.

A partition is a {\it simultaneous} $(s,t)$-core partition if it avoids hook-lengths of both $s$ and $t$. For example, the
above   partition, $(5,4,2,1,1)$, is a $(5,8)$-core partition (and a $(5,10)$-core partition, and a $(100,103)$-core partition etc.)

For a lucid and engaging account, see [AHJ].

As mentioned in [AHJ], Jaclyn Anderson ([An]) very elegantly proved the following.

{\bf Theorem ([An])} If $s$ and $t$ are relatively prime positive integers, then there are
{\bf exactly}
$$
\frac{(s+t-1)!}{s!t!} \quad
$$
$(s,t)$-core partitions. 

For example, here are the $(3+5-1)!/(3!5!)=7$  $(3,5)$-core partitions:
$$
\{ empty , 1, 2, 11, 31, 211, 4211 \} \quad .
$$

{\bf The Order Ideal $P_{n+1,n+2}$}

It turns out that it is most convenient to work with the order ideal 
$$
P_{s,t} := {\bf N} \backslash (s{\bf N} + t{\bf N} ) \quad,
$$
where ${\bf N}$ is the set of non-negative integers. Anderson showed that $(s,t)$-core partitions are in one-to-one correspondence with {\it order ideals} of $P_{s,t}$ ([An]).

Our poset of interest, $P_{n+1,n+2}$, can be identified with a  triangular region in the
2D rectangular lattice, let's call it $A_n$,
$$
A_n:=\{ (i,j) \in {\bf N}^2 \, \vert \, 0 \leq i \leq n-1 \, , \, 0 \leq j \leq n-1-i \} \quad,
$$
consisting of $(n+1)n/2$ lattice points. We {\bf label} the lattice point with
$$
L(i,j) := (n+1)n-1 - (n+2)i - (n+1)j \quad,
$$
in other words, we identify the lattice point $(i,j)$ with the member
 $(n+1)n-1 - (n+2)i - (n+1)j$ of  $P_{n+1,n+2}$.

To see the lattice $P_{9,10}$ see

{\tt http://sites.math.rutgers.edu/\~{}zeilberg/tokhniot/PictOddArmin/O0.html} \quad ,

and to see the lattice $P_{10,11}$ see

{\tt http://sites.math.rutgers.edu/\~{}zeilberg/tokhniot/PictOddArmin/O1.html} \quad .

Note that the point $(n-1,0)$ is labeled $1$, and when we read the labels along diagonals, from the bottom-right to the top-left,
the labels increase by $1$, but as we move from the end of one diagonal to the next one there are ``discontinuities''
of sizes $3,4, \dots,n+1$ respectively.

A subset $S$ of $A_n$ is an {\it order ideal} if it satisfies the following condition:

$\bullet$ If $(i,j) \in S$ then $(i',j') \in S$ for {\it all} $(i',j') \in S$ 
such that $i' \geq i $ and $j' \geq j$. 

In other words, if a lattice point belongs to $S$, then so do all
the lattice points of $A_n$ that are {\it both} to its (weak) right {\it and} are (weakly) above it.
Conversely, and equivalently, 

$\bullet$ If $(i,j) \not \in S$, then the set $\{(i',j') \in A_n \, | \,  i' \leq i \quad and \quad  j' \leq j \}$ is disjoint from $S$.

In other words, if a lattice point is unoccupied, then all the lattice points (weakly) below  and (weakly) to its left
are also unoccupied.

With this interpretation, it is very easy to prove the special case of Anderson's theorem that the
number of $(n+1,n+2)$-core partitions is Catalan's number $(2n+2)!/((n+1)!(n+2)!)$. Given such an order ideal, 
let's call it $S$, let $i$ ($0 \leq i \leq n$) be the smallest positive integer with the property that
$(n-1-i,i)$  is {\bf not} a member of $S$: in other words, the smallest integer $i$ such that
$(n-1,0),(n-2,1), \dots, (n-i,i-1)$ are members of $S$ while $(n-1-i,i)$ is {\bf not} a member of $S$.

Then all the points to the left of and below  $(n-1-i,i)$  are definitely unoccupied, and the order ideal has
two parts, 

$\bullet$ Those (strictly) below and (strictly) to the right of $(n-1-i,i)$  \quad ,

$\bullet$  Those (strictly) above and (strictly) to the left of $(n-1-i,i)$  \quad .

The former component is isomorphic to an order ideal of $A_{i-1}$ with its outer diagonal
fully occupied; i.e., by definition of $i$, $\{(n-1,0), (n-2,1), \dots, (n-i,i-1)\}$ are occupied,
and removing these ``mandatory" members, the remaining set is isomorphic to an order ideal of $A_{i-2}$. 
Let's call it $S_1$.

The top part is isomorphic to an order ideal of $A_{n-i}$, let's call it $S_2$.

This introduces a {\bf canonical decomposition}
$$
S \rightarrow (i, S_1, S_2) \quad, \quad 0 \leq i \leq n \quad , \quad S_1 \in A_{i-2} \quad  , \quad S_2 \in A_{n-i} \quad 
\eqno{(Canonical Decomposition)}
$$
that is obviously one-to-one.
Let $a_n$ be the number of order ideals of $A_n$; then it follows that
it satisfies the recurrence
$$
a_n \, = \, \sum_{i=0}^{n} a_{i-1} a_{n-i-1} \quad,
$$
with the initial condition $a_{-1}=1$. As is well-known (and easy to see) this implies that
indeed $a_n=(2n+2)!/((n+1)!(n+2)!)$.

Note that for the above argument the labels are irrelevant.

{\bf Counting subfamilies}

{\bf Distinct Partitions}

What about counting $(n+1,n+2)$-core partitions into {\bf distinct} parts?
It was conjectured by Amdeberhan[Am], and first proved (as a special case of a much  more general result)
by Straub[S1] that this number is $F_{n+2}$, where $(F_n)$ is the Fibonacci sequence.. Using order ideals, this is immediate. 
By the transformation
$$
(a_1, a_2, \dots , a_k) \rightarrow (a_1-(k-1), \dots, a_k) \quad,
$$
from the sorted list of labels to an $(n+1,n+2)$-core partition, the condition of being distinct translates to the fact
that the corresponding order ideal can't have any adjacent points,
when ``read" along diagonals, from right-to-left and from bottom-to-top. 
This precludes any member of an inner diagonal
(since their existence would imply at least two adjacent points in the outermost diagonal), and the
members of $S$ that do belong to the outer-diagonal can't be adjacent. Let $d_n$ be their number.
If $(n-1,0)$ is not a member of $S$ then $S$ can be viewed as an order ideal (with the above condition)
of $A_{n-1}$, accounting for $d_{n-1}$ such creatures. On the other hand, if $(n-1,0) \in S$, then
$(n-2,1) \not \in S $, and removing both of these yields an order ideal of $A_{n-2}$ (with the above conditions),
hence $d_n$ satisfies the recurrence
$$
d_n = d_{n-1} + d_{n-2} \quad ,
$$
subject to the initial conditions $d_{0}=1$, $d_{1}=2$.

For multi-cores see the elegant paper [AmL].

{\bf Intermezzo: The Joy and Agony of Enumerative Combinatorics}

It is both fascinating and frustrating that in enumeration problems, `tweaking' a problem ever so slightly
turns it from almost trivial (and often, utterly trivial) to very difficult (and often, intractable).
For example, it is utterly trivial that the number of $n$-step walks in the 2D rectangular lattice is
$4^n$, but just add the adjective ``self-avoiding"---in other words, the number of such walks that never
visit the same vertex twice---and the enumeration problem becomes (most probably) intractable and,
at any rate, wide open.

Another example is counting permutations that avoid a pattern. 
The number of permutations, $\pi$, of length $n$ that avoid the pattern $12$ 
(i.e. you can't have $1 \leq i_1<i_2 \leq n$ such $\pi_{i_1}<\pi_{i_2}$) is trivially $1$.
A bit less trivially, but still very doable, is the fact that 
the number of permutations, $\pi$, of length $n$ that avoid the pattern $123$ 
(i.e. you can't have $1 \leq i_1<i_2 <i_3\leq n$ such $\pi_{i_1}<\pi_{i_2}<\pi_{i_3}$) is the good old
Catalan number $(2n)!/(n!(n+1)!)$. But for most patterns, such an enumeration is (probably) intractable.
The simplest wide open case, that we believe is intractable (but we would be happy to
be proven wrong) is to count permutations that avoid the pattern $1324$
(OEIS sequence A061552 [{\tt https://oeis.org/A061552}]), for which the current record is knowing the $36$ first terms.

Returning to  the main topic, consider  enumerating $(2n-1,2n+1)$-core partitions into distinct parts.
Armin Straub conjectured the deceptively simple formula $4^n$.
Alas, its (known) proofs are far from simple!
Straub's conjecture was first proved, 
by  Sherry H.F. Yan, Guizhi Qin, Zemin Jin, Robin D.P. Zhou [YQJZ],
via an ingenious but rather complicated combinatorial proof.
A shorter, but still non-trivial proof was given in [ZZ], using ``guess-and-check," and this was further simplified
by Straub (see [ZZ]). As far as we know, enumerating $(s,t)$-core partitions into distinct parts
for other cases, say $(3n-1,3n+1)$-core partitions, is wide open.

{\bf The problem of counting $(n+1,n+2)$-core partitions into odd parts}

Leonhard Euler famously proved that the number of partitions of an integer $n$ into {\bf distinct} parts equals
the number of partitions of the same $n$ into {\bf odd parts}. (This classical theorem was
recently {\it refined} in a new, very surprising way, by Armin Straub ([S1])). 

Moving on to counting $(n+1,n+2)$-core partitions into {\it odd} parts 
(note that we are counting {\it all} of them, regardless of their size, so the $n$ here is not the same as the $n$
in Euler's theorem), it seems that the number of such partitions has nothing to do with
the number of $(n+1,n+2)$-core partitions into  {\it distinct} parts (i.e. $F_{n+2}$).
This new problem seems (at least to us) much harder. On the other hand, we believe that it should
be {\it doable}, and invite anybody to  tackle it (see the abstract) and thereby support the OEIS.

We will now describe our approach, its success (it enabled us to crank out $23$ terms, thereby extending
Straub's $11$ terms, and with better computers, and more optimization, one may be able to crank out a few more terms),
and its major shortcoming. At the end of the day, it is an {\it exponential time} (and memory!) algorithm.

{\bf A Scheme for Counting Order Ideals of $P_{n+1,n+2}$ Corresponding to $(n+1,n+2)$-core partitions into odd parts}

It is readily seen, by the mapping from order ideals to partitions
$$
(a_1, \dots , a_k) \rightarrow (a_1+k-1, a_2+k-2, \dots, a_k+0) \quad
$$
that  an order ideal of $P_{n+1,n+2}$ corresponds to an $(n+1,n+2)$-core partition into {\bf odd} parts
if and only if, when reading the occupied labels, as described above, along diagonals, from bottom-right to top-left, starting
from the rightmost diagonal and ``walking" to the left, (i) the first label read is odd and (ii) the labels alternate in parity. Since only
the parity matters, we can color the vertices of $A_n$  by the colors ``even" and ``odd."

Alas, one has to distinguish two cases.  For both $n$ even and odd, the label of $(n-1,0)$ is odd
(since it is always $1$), and as you
proceed, in $A_n$ along diagonals, the parities {\bf alternate}. But for $n$ odd, all the parities
along the same row are the same, while if $n$ is even, they alternate. Hence we are forced to consider
the more general problem where  there is a ``coloring" parameter, let's call it $c$, ($c=0$ or $c=1$)
such that the ``color" is
$$
C(i,j) := 1 + \, c\,i \, + \, (1-c)j \quad modulo \quad 2 \quad.
$$

So let's forget, for now, about $(n+1,n+2)$-core partitions into odd parts, and  instead define the following:

$\bullet$ Let $e^{(0)}(n)$ be the number of order ideals of $P_{n+1,n+2}$ such that when read along diagonals,
the occupied vertices alternate in color using coloring parameter $c=0$, and the first label is odd.

$\bullet$ Let $e^{(1)}(n)$ be the number of order ideals of $P_{n+1,n+2}$ such that when read along diagonals,
the occupied vertices alternate in color using coloring parameter $c=1$, and the first label is odd.

Once we would find  a way to compute both sequences $e^{(0)}_n$ and $e^{(1)}_n$,
our object of desire, the Straub sequence, enumerating $(n+1,n+2)$-core partitions into odd parts, let's call
it $s_n$, would be given by 
$$
s_n \, = \, \cases{ e^{(0)}_n \, , \, if \quad n \quad is \quad even \, ; \cr
                    e^{(1)}_n \, , \, if \quad n \quad is \quad odd \, } \quad .
$$

For two examples of order ideals corresponding to $(10,11)$-core partitions into odd parts, see

{\tt http://sites.math.rutgers.edu/\~{}zeilberg/tokhniot/PictOddArmin/O2.html} \quad ,

and

{\tt http://sites.math.rutgers.edu/\~{}zeilberg/tokhniot/PictOddArmin/O3.html} \quad .

{\bf Dynamical Programming}

The obvious way would be to try and extend the  above argument for counting {\it all} order ideals of $P_{n+1,n+2}$,
using $(CanonicalDecomposition)$.

Let  $(n-1-i,i)$ be the `first' unoccupied vertex of the order ideal $S$ of $A_n$,
let $(i,S_1,S_2)$ be its image under $(CanonicalDecomposition)$.
The smaller order ideals $S_1$ (of $A_{i-2}$) and $S_2$ (of $A_{n-i}$)
also have the property, that 
{\it within} each diagonal, the colors of the occupied vertices alternate, but, alas, as you move
from one diagonal to the next one, the ``alternation" may (and often does) ``break-down". Also the
two components in the ``canonical decomposition" are not ``independent" but must satisfy some
compatibility conditions.

This forces us to consider much more general creatures, order ideals whose ``colors" (parity) alternate
within each individual diagonal, and having, {\bf in addition}, a given ``coloring profile", the 
list of pairs of colors of the first and
last vertices in each diagonal, reading from left to right.  For example, the profile of the order ideal in

{\tt http://sites.math.rutgers.edu/\~{}zeilberg/tokhniot/PictOddArmin/O2.html}

is $[[1,1],[0,0],[1,0]]$, since it is supported in the three outermost diagonals and the 
occupied vertices on rightmost diagonal start with $3$ (odd parity, hence $1$), and end with $9$
(hence the first component of the profile is $[1,1]$). The lowest occupied vertex on the second
diagonal has label $14$ (hence $0$) and the last one has label $18$ (hence $0$), hence for the
second diagonal, we have $[0,0]$. Finally, the lowest label on the third diagonal is $25$
and the highest is $26$, hence $[1,0]$. Note that for any profile of an order ideal corresponding
to an $(n+1,n+2)$-core partition
$$
[[a_1,b_1], [a_2,b_2], \dots , [a_k,b_k]] \quad ,
$$
$b_i$ and $a_{i+1}$ must have {\bf opposite} parities. Also, $a_1=1$.  We call such profiles {\bf good profiles}. Hence there are $2^{k-1}$ good profiles. Unfortunately, in order to use
dynamical programming, we need to consider all $2^{2k}$ profiles for $k$ diagonals
(and it is easy to see that for us, $k \leq n/2$). Hence our algorithm is {\it exponential} time
(and memory).

We essentially use $(Canonical Decomposition)$ but refined to order ideals with a given profile,
and at the end we sum over all good profiles.

The details are straightforward but rather tedious, and may be 
gotten from looking at the source code of the Maple package {\tt OddArmin.txt}, available from

{\tt http://www.math.rutgers.edu/\~{}zeilberg/tokhniot/OddArmin.txt} \quad .

See procedure {\tt NuOIG(n,c)} giving $e^{(c)}_n$ for $c=0$ and $c=1$. It is obtained by adding up
the outputs of procedure {\tt NuOIP(n,c,P)} where {\tt P} is the `profile', and adding up all
the outputs from the set of `good profiles'. 

Procedure {\tt NuOIP(n,c,P)} is a refined version of $(CanonicalDecomposition)$.
The problem is the {\it proliferation} of profiles. For each ``good" order ideal with  given profile $P$, for which 
$(n-1-i,i)$ is the `first' unoccupied vertex, the two smaller order ideals, of  $A_{i-2}$ and $A_{n-i}$, have
implied profiles. The program finds all such pairs {\tt (P1,P2)} and adds up 

{\tt NuOIP(i-2,c,P1)*NuOIP(n-i,c,P2)} \quad,

for all such pairs.

The output was as follows.

$\bullet$ The first $23$ terms of the sequence $e^{(0)}_n$ (staring with $n=0$) are
$$
2, 4, 7, 17, 30, 80, 143, 404, 728, 2140, 3876, 11729, 21318, 65952, 120175,
$$
$$
    378321, 690690, 2205168, 4032015, 13023324, 23841480, 77761008, 142498692 \quad .
$$

$\bullet$ The first $23$ terms of the sequence $e^{(1)}_n$ (staring with $n=0$) are

$$
2, 3, 7, 12, 31, 55, 152, 273, 790, 1428, 4271, 7752, 23767, 43263, 135221, 46675, 782968, 
$$
$$
1430715, 4598804, 8414640, 27332956, 50067108, 164081764
$$

But, we really don't care about $e^{(0)}_n$ when $n$ is odd, or $e^{(1)}_n$ when $n$ is even. We want 
the Straub sequence $e^{(n \, mod \,\,2)}_n$. In other words, we extract the even-indexed terms of the former sequence
and the odd-indexed terms of the latter sequence, and then we {\it interleave} them.
This yields the first $23$ terms of the Straub sequence:
$$
1,2, 4, 7, 17, 31, 80, 152, 404, 790, 2140, 4271, 11729, 23767, 65952, 135221, 
$$
$$
378321, 782968, 2205168, 4598804, 13023324, 27332956, 77761008 \quad .
$$

This sequence is {\bf not} (yet) in the OEIS (but hopefully it will be very soon). But what about the ``rejected" terms,
the ones that we do not care about? Maybe we {\bf should} care about them!

The  first $23$ terms of the sequence  $e^{(n+1 \, mod \,\,2)}_n$ are
$$
1, 2, 3, 7, 12, 30, 55, 143, 273, 728, 1428, 3876, 7752, 21318, 43263, 120175, 
$$
$$
246675, 690690, 1430715, 4032015, 8414640, 23841480, 50067108 \quad .
$$

To our utter surprise (and delight), this sequence {\bf is} in the OEIS, (but for entirely different reasons!).  It is 
sequence A047749 and has a very nice closed-form expression: If $n=2m$, then $\frac{1}{2m+1} \cdot {{3m} \choose {m}}$, while if $n=2m+1$, then $\frac{1}{2m+1} \cdot {{3m+1} \choose {m+1}}$.  As mentioned in the OEIS entry, it is easily verified that its generating function, $Y=Y(x)$, satisfies the simple cubic equation
$$
x\, Y^3 \,- \, 2\,Y^2 \,+ \, 3\,Y \,-\, 1 \,= \, 0 \quad .
$$

We are almost sure that the generating function of the Straub sequence $s_n=e^{(n \, mod \,\,2)}_n$ also satisfies an algebraic equation,
but the above $23$ terms did not allow us to guess one.

As we mentioned in the abstract, we would gladly donate one hundred dollars to the OEIS, in honor of the first person
to generate enough terms that would enable the discovery of such an algebraic equation (with a few terms to spare, yielding
a non-rigorous proof), and an additional one hundred dollars (either in honor of the same or  different person(s) and/or machines),
for a rigorous proof.

{\bf Enumerating Restricted Families of Core Partitions}

Finally, we shall discuss some bonus families of partitions related to Straub's paper.  In these cases, we were able to use symbolic computation to rigorously derive rational generating functions.

{\bf $(n+1,n+2)$-core partitions with at most $k$ repeats of a part}

As noted, the sequence of numbers enumerating 
$(n+1,n+2)$-core partitions with distinct parts is $\{ F_{n+2}\}_{n=0}^{\infty}$, whose
generating function is the very simple rational function $\frac {1+x}{1-x-{x}^{2}}$.
We shall now show that it is not hard to derive such rational generating functions to enumerate
$(n+1,n+2)$-core partitions where  each part gets repeated at most $k$ times, for any, given, {\it specific} (i.e. numeric, not symbolic)
$k$, where the former case corresponds to $k=1$.

Again, consider the poset $A_n:=P_{n+1,n+2}$, whose order ideals correspond to $(n+1,n+2)$-core partitions. Suppose $S$ is an order ideal of $A_n$ corresponding to a partition in which each part appears at most $k$ times. This is equivalent to saying $S$ contains at most $k$ consecutive labels. (Note that, because $S$ is an order ideal, a {\it necessary} condition for this is that the elements of $S$ reside in the $k$ outermost diagonals of $A_n$.)

As before, let $(n-1-i,i)$ be the smallest-labeled unoccupied point in the outermost diagonal of $S$, so that $S$ contains the labels $1,\dots,i$ but not $i+1$. Due to our new restriction, $i\leq k$. Again, let $S_1$ contain the elements of $S$ below $(n-1-i,i)$ and not on the outer diagonal; let $S_2$ contain elements above $(n-1-i,i)$. Then $S_1$ is isomorphic to an arbitrary order ideal of $A_{i-2}$, and $S_2$ is isomorphic to an order ideal of $A_{n-i}$ with no $k$ consecutive labels.

So, with $k$ fixed, we can see $S$ as the ``concatenation'' of two types of order ideals---one with a filled-in base of size $\leq k$, and another of the same type as $S$. The generating function enumerating the first type of order ideals is a finite polynomial: its coefficients are Catalan numbers. So we obtain an algebraic equation satisfied by the desired generating function that can easily be solved in Maple. See the procedure {\tt Fk} in the Maple package.

Here are the generating functions for $2 \leq k \leq 4$:

$k=2$: 
$$
-{\frac {2\,{x}^{2}+x+1}{2\,{x}^{3}+{x}^{2}+x-1}} \quad ,
$$
whose first few coefficients are
$$
1, 2, 5, 9, 18, 37, 73, 146, 293, 585, 1170, 2341, 4681, 9362, 18725, 37449,\dots, \quad ;
$$

$k=3$:
$$
-{\frac {5\,{x}^{3}+2\,{x}^{2}+x+1}{5\,{x}^{4}+2\,{x}^{3}+{x}^{2}+x-1}} \quad ,
$$
whose first few coefficients are
$$
1, 2, 5, 14, 28, 62, 143, 331, 738, 1665, 3780, 8576, 19376, 43837, 99265, 224734,\dots \quad ;
$$

$k=4$:

$$
-{\frac {14\,{x}^{4}+5\,{x}^{3}+2\,{x}^{2}+x+1}{14\,{x}^{5}+5\,{x}^{4}+2\,{x}^{3}+{x}^{2}+x-1}} \quad ,
$$
whose first few coefficients are
$$
1, 2, 5, 14, 42, 90, 213, 527, 1326, 3317, 8022, 19608, 48272, 119073, 293109,
$$
$$
719074, 1766201, 4342666, 10679582, 26253546, 64516501, 158569355, 389788182 \dots \quad .
$$

For the generating functions, and first few terms, for the cases $5 \leq k \leq 20$, see the webpage

{\tt http://sites.math.rutgers.edu/\~{}zeilberg/tokhniot/oOddArmin3.txt} \quad .

{\bf $(n+1,n+2)$-core partitions into odd parts, whose corresponding order ideals are confined to the $k$ outermost diagonals}

Inspired by the necessary condition mentioned above, let us enumerate $(n+1,n+2)$-core partitions into {\it odd} parts whose order ideals are restricted to the outer $k$ diagonals.

As before, we classify $S$ according to its profile $P$, a list of pairs, each pair giving the parities of the largest and smallest labels of $S$ in a certain diagonal. Also, define $i(S)$ to be the smallest $j$ such that $(j,0)$ is occupied. 

Call $(P(S),i(S))$ the ``type'' of $S$; for fixed $k$, there are finitely many types. Further, any $S$ of a certain type is the concatenation of its elements on the $x$-axis with some smaller order ideal of $A_{n-1}$ having a compatible type. Thus, the generating function of order ideals having a certain type satisfies some algebraic equation involving the generating functions of its ``child'' types. Once we solve this system and sum the generating functions over $P$, we get what we are after. See {\tt Gk} in the Maple package.

For $k=2$ the generating function is
$$
-{\frac {{x}^{4}-{x}^{3}-{x}^{2}+x+1}{{x}^{5}-{x}^{4}-2\,{x}^{3}+3\,{x}^{2}+x-1}} \quad,
$$
whose first few coefficients are
$$
1, 2, 4, 7, 15, 27, 56, 104, 210, 398, 791, 1517, 2988, 5769, 11306, 21911, 
$$
$$
42820, 83160, 162261, 315496, 615050, 1196676, 2331733, 4538426, 8840719 \quad  .
$$

For $k=3$ the generating function is
$$
-{\frac {{x}^{9}+{x}^{8}-4\,{x}^{7}-6\,{x}^{6}+8\,{x}^{5}+9\,{x}^{4}-5\,{x}^{3}-5\,{x}^{2}+x+1}{ \left( {x}^{9}+2\,{x}^{8}-3\,{x}^{7}-9\,{x}^{6}+3\,{x}^{5}+
14\,{x}^{4}-{x}^{3}-7\,{x}^{2}+1 \right)  \left( x-1 \right) }} \quad ,
$$
whose first few coefficients are
$$
1, 2, 4, 7, 17, 31, 76, 144, 344, 670, 1560, 3103, 7079, 14315, 32152, 65861,
$$
$$
146183, 302456, 665300, 1387172, 3030464, 6356068, 13813464, 29103412, 62999146
\dots.
$$

For generating functions, and first few terms, for the cases $4 \leq k \leq 5$, see the webpage

{\tt http://sites.math.rutgers.edu/\~{}zeilberg/tokhniot/oOddArmin2.txt} \quad .

{\bf References}

[Am] Tewodros Amdeberhan, {\it Theorems, problems and conjectures}, \hfill\break
{\tt https://arxiv.org/abs/1207.4045} \quad .

[AmL] Tewodros Amdeberhan and  Emily Sergel Leven , {\it Multi-cores, posets, and lattice paths}, \hfill\break
{\tt https://arxiv.org/abs/1406.2250}.
Also published in Adv. Appl. Math. {\bf 71}(2015), 1-13.

[An] Jaclyn Anderson, {\it Partitions which are simultaneously $t_1$  and $t_2$ -core}, Discrete Math. {\bf 248}(2002), 237-243.

[AHJ] Drew Armstrong, Christopher R.H. Hanusa, and B. Jones, {\it Results and conjectures on simultaneous core partitions}, 
{\tt https://arxiv.org/abs/1308.0572} \quad .
Also published in European J. Combin. {\bf 41} (2014), 205-220.

[J] Paul Johnson, {\it Simultaneous cores with restrictions and a question of Zaleski and Zeilberger}, \hfill\break
{\tt https://arxiv.org/abs/1802.09621} \quad .

[S1] Armin Straub, {\it Core partitions into distinct parts and an analog of Euler's theorem}, \hfill\break
{\tt https://arxiv.org/abs/1601.07161}.
Also published in: European J. of Combinatorics {\bf 57} (2016), 40-49.

[S2] Armin Straub, {\it  Core partitions into distinct parts and an analog of Euler's theorem (JMM)}, (talk),
{\tt http://arminstraub.com/talk/corepartitions-jmm} \quad .

[YQJZ] Sherry H.F. Yan, Guizhi Qin, Zemin Jin, Robin D.P. Zhou, {\it On $(2k+1,2k+3)$-core partitions with distinct parts}, 
{\tt https://arxiv.org/abs/1604.03729} \quad .

[Za] Anthony Zaleski, {\it Explicit expressions for the moments of the size of an $(s,s+1)$-core partition with distinct parts},
{\tt https://arxiv.org/abs/1608.02262}. 
Also appeared in Advances in Applied Mathematics Vol {\bf 84} (2017), 1-7.

[ZZ] Anthony Zaleski and Doron Zeilberger, 
{\it Explicit (Polynomial!) Expressions for the Expectation, Variance and Higher Moments of the Size of of a (2n+1,2n+3)-core partition with Distinct parts  },
{\tt https://arxiv.org/abs/1611.05775} and \hfill\break
{\tt http://sites.math.rutgers.edu/\~{}zeilberg/mamarim/mamarimhtml/armin.html} \quad .
Also appeared in J. Difference Equations and Applications Volume {\bf 23}(2017), 1241-1254, \hfill\break

\bigskip
\hrule
\bigskip
Anthony Zaleski, Department of Mathematics, Rutgers University (New Brunswick), Hill Center-Busch Campus, 110 Frelinghuysen
Rd., Piscataway, NJ 08854-8019, USA. \hfill\break
Email: {\tt az202 at math dot rutgers dot edu}   \quad .
\bigskip
Doron Zeilberger, Department of Mathematics, Rutgers University (New Brunswick), Hill Center-Busch Campus, 110 Frelinghuysen
Rd., Piscataway, NJ 08854-8019, USA. \hfill\break
Email: {\tt DoronZeil at gmail dot com}   \quad .
\bigskip
Exclusively Published in the Personal Journal of Shalosh B. Ekhad and Doron Zeilberger and arxiv.org  \quad .
\bigskip
First version: Dec. 26, 2017; Previous version: Jan. 24, 2018; This version: Feb. 28, 2018.

\end